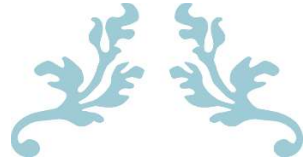

# INTERNATIONAL JOURNAL OF MECHANICAL AND INDUSTRIAL TECHNOLOGY

## ISSN 2348-7593

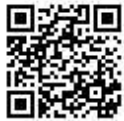

https://www
.researchpublish
.com/journal-detail

# Innovative Line Balancing for the Aluminium Melting Process

Prof. Dr. Ray Wai Man Kong[1], Ding Ning[2], Theodore Ho Tin Kong[3]

[1]Adjunct Professor, System Engineering Department, City University of Hong Kong, China

[1]Eagle Nice International (Holding) Ltd., Hong Kong, China

[2]Engineering Doctorate Student, System Engineering Department, City University of Hong Kong, China

[3] Graduated Student, Master of Science in Aeronautical Engineering, Hong Kong University of Science and Technology, Hong Kong

[3] Thermal-acoustic (Mechanical) Design Engineer at Intel Corporation in Toronto, Canada

*Abstract:*  This research article explores the optimization of aluminium extrusion processes through advanced line balancing techniques, focusing on maximizing marginal profit by increasing melting and casting outputs. By employing mixed integer linear programming (MILP), we identify strategies to minimize idle costs and enhance production efficiency. The study demonstrates that increasing the daily cycle rate from 2 to 4.36 cycles results in a significant rise in daily marginal profit, calculated at USD67,786, after accounting for additional labour costs. This optimization is achieved by expanding the workforce from 8 to 12 operators across two shifts, leading to a 50% increase in labour expenses. The findings reveal a remarkable 117.6% growth in marginal daily profit, underscoring the potential of automation and intelligent manufacturing in transforming the aluminium extrusion industry. Insights from cross-industry research, including Lean Methodology in the Modern Garment Industry, further illustrate the broader applicability of these advancements. This study highlights the critical role of automation in driving productivity and profitability in manufacturing sectors, paving the way for future innovations in aluminium extrusion and beyond.

## I.  INTRODUCTION

As a global market research company, the Aluminium Extrusion Market Size was valued at USD 83.9 Billion in 2023 from Market Research Future (MRFR) [1] in Fig. 1, which is from the MRFR Database and Analyst Review. The Aluminium Extrusion industry is projected to grow from USD 90.77 Billion in 2024 to USD 170.53 Billion by 2032, exhibiting a compound annual growth rate (CAGR) of 8.20% during the forecast period (2024 - 2032). The Aluminium Extrusion Market will be driven by a rise in the need for strong, lightweight extruded items with high corrosion resistance, and economic expansion, accelerated urbanisation and expanding infrastructure projects are the key market drivers enhancing market growth.

The aluminium extrusion market has been driven by the increasing demand for lightweight and durable extruded products. Aluminium extrusions are widely used in various industries, including construction, automotive, aerospace, electronics, and consumer goods, among others. One of the main advantages of aluminium extrusions is their lightweight nature. Aluminium has a high strength-to-weight ratio, making it an ideal material for applications where weight reduction is a critical factor, such as in the automotive and aerospace industries. Aluminium is highly durable and corrosion-resistant, which makes it an attractive choice for products that need to withstand harsh environments or extreme weather conditions. The aluminium frame can be used for the construction of the required material for the iPad, iPhone and window. The major materials of the auto vehicles are the aluminium alloy frame.

The aluminium alloy manufacturer has the most critical role in the extrusion process, which allows for the creation of complex shapes and designs, which can be used to meet the specific needs of different industries and applications.

## II.   MANUFACTURING PROCESS OF ALUMINIUM ALLOY

*Manufacturing Process of Melting for Aluminium Alloy Manufacture*

Aluminium Alloy Melting is a major process in which a solid heats up and becomes a liquid state of aluminium. After the feeding is completed, melting starts. During the smelting process, it is necessary to ensure that the melting is fast and uniform and that the temperature near the flame is high, reaching more than 1200 °C. The melting temperature is reaching over 700°C. To avoid local overheating and excessive temperature, the aluminium alloy oxidation causes a serious inconvenience to the later product refining. Stirring should be carried out, and the unmelted furnace charge should be squeezed into the aluminium liquid so that almost all of it is immersed in the stirring to avoid local overheating. For the final immersion of Magnesium (Mg) and other metallic elements for various aluminium alloys, the flame cannot be directly heated and melted, and because the aluminium liquid is immersed in the raw material, the temperature is reduced, and the Mg and other metallic elements are melted at a relatively low temperature, which reduces the burning loss and improves the combustion efficiency.

The melting process chart in Fig. 2 is shown below to relate to the melting process. The machinery standard time and labour standard time are the critical time measurements for the line balancing of the manufacturing process. The input substance is the raw aluminium and recycled aluminium, and the output is the aluminium rod. The aluminium rod is a large and heavy state of aluminium to direct extrusion for reforming the various shapes of aluminium to the related manufacturing processes for the automobile, the frame of iPad & iPhone and window frame.

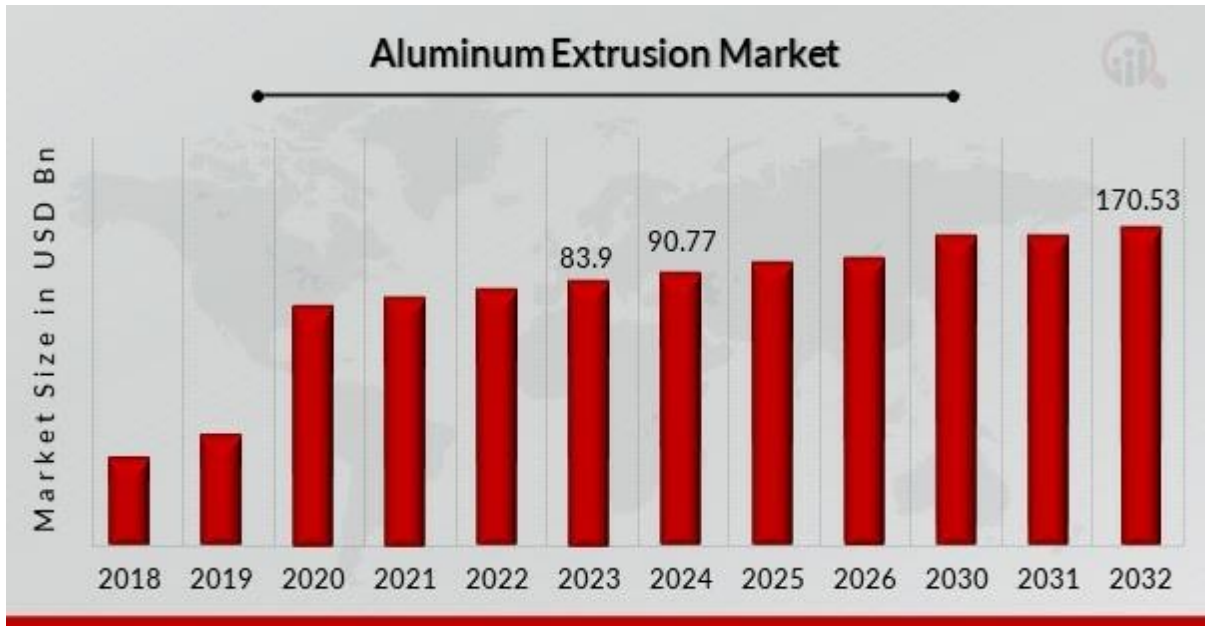

**Figure 1: Aluminium Extrusion Market Research from MRFR Database and Analyst Review**

*The problem of Line Unbalancing in Aluminium Melting*

Line balancing in the aluminium manufacturing industry is the technique of levelling the output of every operation in an aluminium alloy production line, so the output from the upstream operation can be optimized to pass through the downstream operation. There should neither be an accumulation of work between two processes (operation) nor a shortage of workpieces from the upstream workstation (previous work step) between the melting line and its inter-process. It is important to maintain this balance because in a melting process line, the output of one process, as an aluminium rod, is the input of the next aluminium extrusion.

A melting furnace and melting tank for stabilisation are not balanced; hence, there would be the following production problems:

- Reduced Efficiency:

  It means that an upstream aluminium melting output is a downstream operation input. Because of this reason, the worker after the melting process will not get the loading input as per their capacity of producing output, hence they will be underutilized. In this case, it is to make matters worse, more machines and manpower will be allocated to increase aluminium rod sewing and cutting, but efficiency will fall even more.

- Reduced the utilisation and efficiency of the melting furnace:

  The melting furnace has poured the aluminium liquid into the melting tank, and the melting furnace is idle to wait for the melting tank for several processes to the casting.

- Energy wastage for keeping the high temperature in the melting furnace:

  The melting furnace is required to keep the high temperature and stir the residue of the melting furnace. Once the melting furnace is lower than the melting point of aluminium, the aluminium solidifies to a solid state as steady stabilization of the alloy structure. It is not easy to melt and mix the new lot of aluminium or recycled aluminium.

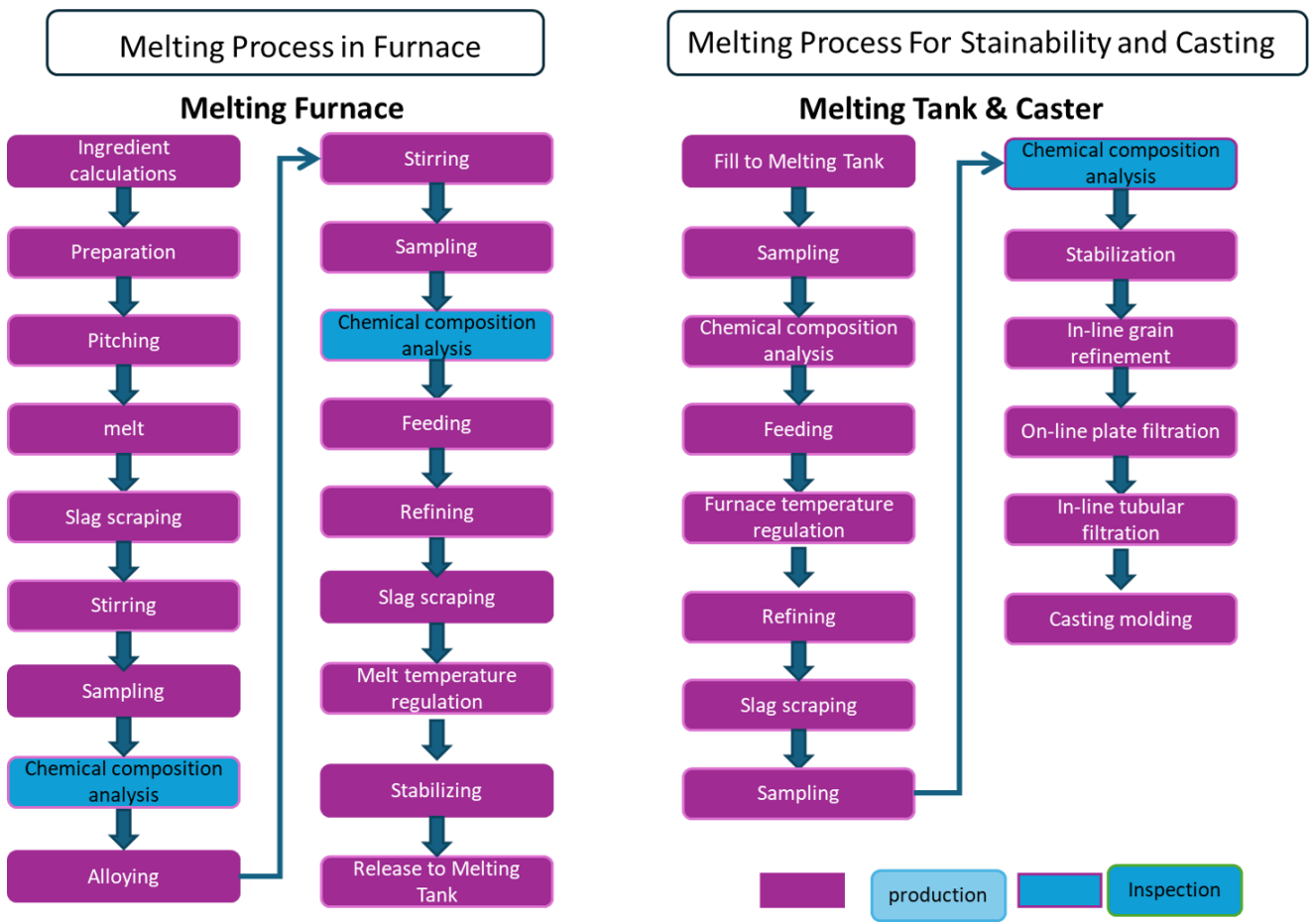

**Figure 2: Melting Process Chart for Aluminium Alloy Manufacturing**

In the melting process, the operations in the melting furnace in Fig. 3 are defined as the melting workstation. The operations in the melting tank are defined as the melting stabilization workstation in Fig. 4. The line balancing is required to make the balance of the melting workstation and the stabilization workstation because these are related to the separated facility and equipment. The workstation of stabilization workstation includes the aluminium rod casting in Fig. 4.

Page | 75

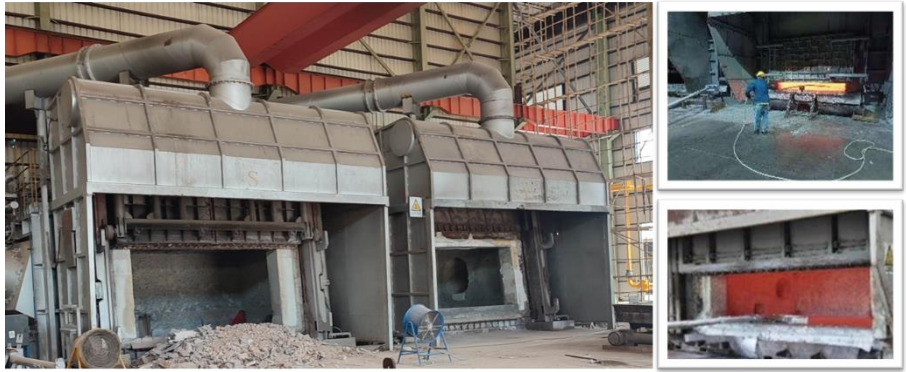

**Figure 3: Melting Furnace for Aluminium Alloy Manufacturing**

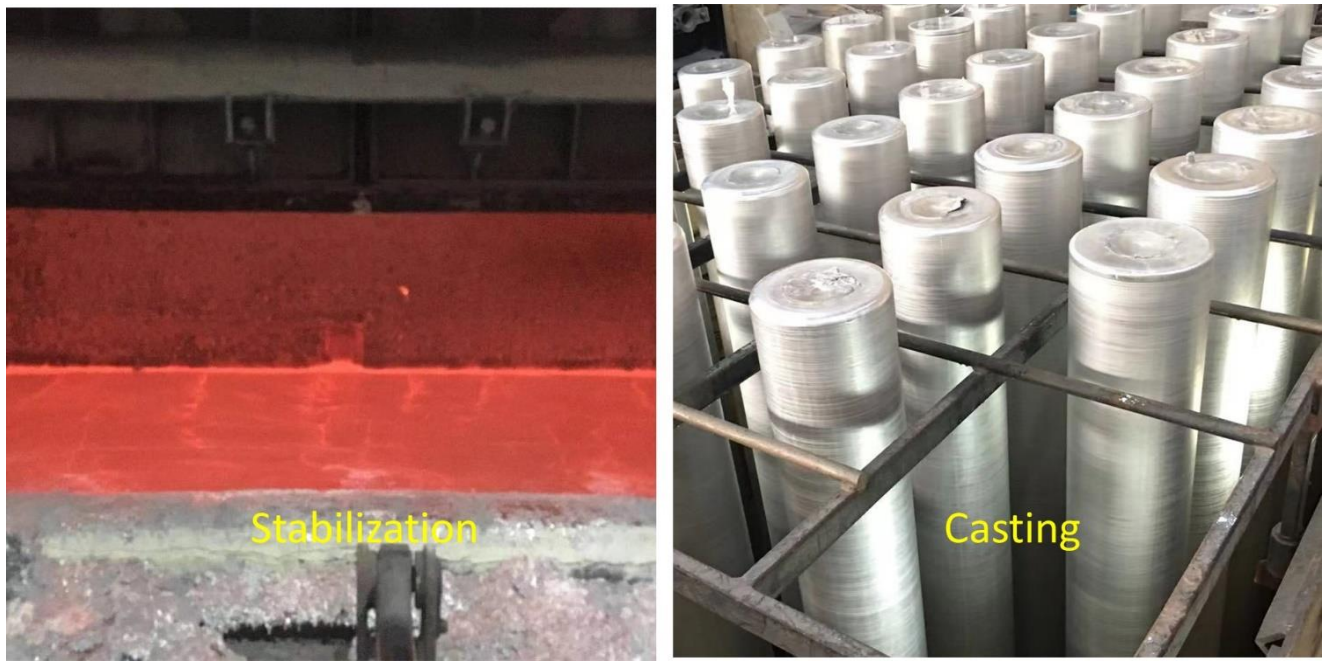

**Figure 4: Melting Stabilization and Casting for Aluminium Alloy Manufacturing**

The traditional setup of the melting furnace line with the melting tank ratio is one-to-one, which is not followed by the line balancing concept. The study of an aluminium extrusion company is proposed for the use of the standard time ratio between the melting furnace and the melting tank. The below section 4 shows the case study how an aluminium company refer to the publication Journals from Prof Dr Ray WM Kong et al [2]. The imbalanced line balance in the melting line affects the shortage of supply of aluminium rods to the aluminium frame machining assembly. For the aluminium assembly, the upstream operation is the supply of aluminium rod and the downstream operation is the aluminium frame machining process. It is not balanced the line balance; hence, the following shortage of supply causes these problems have exist as shown above points.

### III. LITERATURE REVIEW

Prof. Dr Ray Wai Man Kong [3] has discussed strategies to optimize the mixed integer linear programming for line balancing and balance the capacity of machines, machine centres, and work centres in the initial stages of line balancing to enhance output rates. In the context of aluminium extrusion, merely increasing the capacity of individual machines or assembly lines does not necessarily improve overall production output and productivity due to line-balancing challenges. The article "Lean Methodology for Lean Modernization" provides insights into applying lean technology to develop future state value stream mapping (VSM) and identify bottlenecks in the aluminium extrusion process, thereby enhancing capacity and achieving balanced production.

In the article on mixed integer linear programming in the Garment Line Balancing, there is a focus on the mathematical optimization method. For line balancing, Ocident Bongomin [4], the Assembly Line Balancing Problem (ALBP), also known as assembly line design, is a family of combinatorial optimization problems widely studied for its simplicity and industrial applicability. ALBP is NP-hard, encompassing the bin packing problem as a special case. ALBPs arise whenever an assembly line is configured, redesigned, or adjusted. In aluminium extrusion, these problems are particularly pronounced due to the unique balancing challenges posed by the extrusion process compared to other manufacturing lines, such as those for trucks, buses, or machinery.

The task involves distributing the total workload for producing aluminium profiles among workstations along the line, adhering to strict or average cycle times. The general principles of line balancing include considerations for machining, assembly, and disassembly environments; the number of product models (single-model, mixed-model, multimodal lines); and line layouts (basic straight lines, U-shaped lines, circular transfer lines).

The assembly line balancing (ALB) problem has been extensively studied, as noted by Gary Yu-Hsin Chen [5]. The ALB model ensures that staff assignments balance the entire production process, effectively reducing production time or idle time. In aluminium extrusion, the mastery of skills by employees at each task is a critical indicator for achieving ALB. However, there is limited research on multifunctional workers with varying skill levels at workstations. Our research incorporates the Toyota Production System (TPS) principles, adapted for aluminium extrusion, to optimize floor space, flexibility, and working conditions. This approach features U-shaped assembly lines and teams of workers managing extrusion processes on a single-piece flow basis.

Chen et al. [6] address a multi-skill project scheduling problem, which is relevant to aluminium extrusion where projects are divided into tasks completed by skilled employees. Their multi-objective nonlinear mixed integer programming model considers employees' skill proficiency, multifunctional roles, and cell formation to minimize production cycle time. This approach, which accounts for real-world skill variations, effectively reduces production time through optimal personnel assignment and preferred production modes. The human factor introduces uncertainty affecting actual cycle time, emphasizing the need for real-time dynamic line balancing in aluminium extrusion.

Haile Sime & Prabir Jana (2018) [7] demonstrated the use of Arena simulation software to design and evaluate alternative production systems, optimizing resource utilization through effective line balancing. Markus Proster & Lothar Marz (2015) highlighted the importance of dynamic balancing for high productivity in mixed-model assembly lines, applicable to aluminium extrusion where varying extrusion times for different profiles require adaptive strategies. Simulation tools can visualize these methods, reducing complexity and enhancing transparency in planning extrusion lines.

Ghosh and Gagnon (1989), along with Erel and Sarin (1998), provided detailed reviews on these topics. Configurations of extrusion lines for single and multiple products can be divided into single-model, mixed-model, and multi-model types. Single-model lines extrude one product type, mixed-model lines handle multiple products, and multi-model lines produce sequences of batches with intermediate setup operations (Becker & Scholl, 2006).

## IV. METHODOLOGY

*A. Industrial Engineering and Lean Technology to Study the Line Balance of Aluminium Extrusion*

Following the Lean Methodology from the Lean Methodology for the Modern Garment Industry. Industrial engineering applies to the lean methodology and technology in the line balancing of aluminium extrusion manufacturing. An industrial engineer is working for the Here's how it is utilized:

(a) Work Measurement: Industrial engineering study involves conducting time and motion studies to measure the time taken to perform each operation in the melting, stabilisation and casting processes. This cycle time data is related to the machinery for calculating each operation's cycle time and relates to the capacity.

(b) Refining Process in the Melting Analysis: the entire refining process is required to get the analysis report of aluminium alloy composition during the melting process at high temperatures. If the refinement cannot achieve the aluminium alloy requirement, the refinement operation is required to be done again, repeated more times, to achieve the required aluminium alloy composition and strength. The simple process bottlenecks, inefficiencies, and areas of improvement in the refinement process can be identified by Value Stream Mapping. By understanding the process flow, the future state of VSM can identify opportunities for line balancing to optimize the line balancing between the melting furnace operation and the melting tank operation and optimize the refinement time.

(c) Capacity and Manpower Resource Allocation: Industrial engineers assess the aluminium melting and casting workforce and equipment available in the manufacturing facility. They determine the number of operators required for each operation based on the involved machinery cycle time and capacity. The melting process relies on a melting furnace and machinery, but the refining and stirring operations are required for the skilful operators to control the machine and collect the sample of aluminium alloy liquid for the metallic composition and strength test. For the heavy machinery and aluminium industry, it applies more than 1 operator to operate the facility and machinery to stir the aluminium liquid and add the magnesium and other metals to mix the aluminium alloy in the melting furnace. The crows of operators are not fully loaded to operate the furnace and machine. The lean practice with the industrial engineering concept is a good tool to optimize the operators' working time for sharing their time to involve more steps, reducing their waiting time and idle time. This helps in allocating resources effectively and achieving a balanced line.

(d) Layout Design: Industrial engineers consider the layout of the time study of the melting process, stabilization process, refinement process and casting, which has an impact on efficiency and output. They analyze the flow of materials, equipment placement, and operator movement. Optimizing the layout, including batch layout for setting the optimization of the aluminium tank and the casting process. The one melting furnace to one melting tank and casting facility is not an optimization. To ensure the safety of the facility and not solidification of aluminium liquid from the melting furnace to the melting tank, the maximum ratio of the melting furnace to the melting tank is 2 to 1 because of the length of the connection pipe and the liquid flow rate between melting furnace to the melting tank in Fig. 5.

(e) Continuous Improvement: Industrial engineering study emphasizes continuous improvement in line balance as referred to the Mixed Integer linear programming for Garment Line Balancing and Lean Methodology. Industrial engineers monitor the bottleneck operation of the refinery performance, collect data, and analyze it to identify areas for further optimization. They implement changes, conduct follow-up studies, and refine the line-balancing process to achieve higher efficiency and productivity.

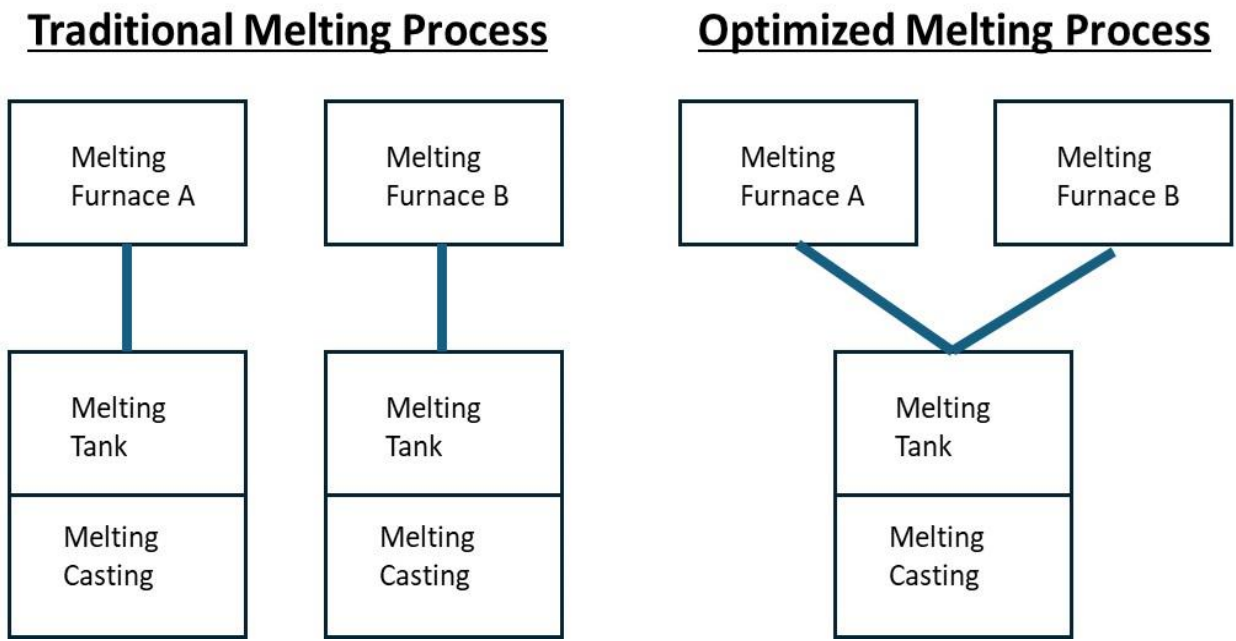

Figure 5: Melting Line Comparison Diagram

## V. CASE STUDY FOR THE LINE BALANCING FOR MELTING PROCESS OPTIMIZATION

By utilizing an industrial engineering study in line balancing aluminium extrusion manufacturing, Company A in the case study can optimize its production processes, reduce lead times, improve resource utilization, and enhance overall efficiency. This results in increased productivity, cost savings, and improved customer satisfaction.

The traditional melting process uses the melting furnace to link with the melting tank and the casting facility in series with a ratio of 1:1. Aluminium liquid in the melting furnace is poured into the melting tank by direct pipe connection. The time study of the traditional method was 2 times a day as shown in Fig. 6.

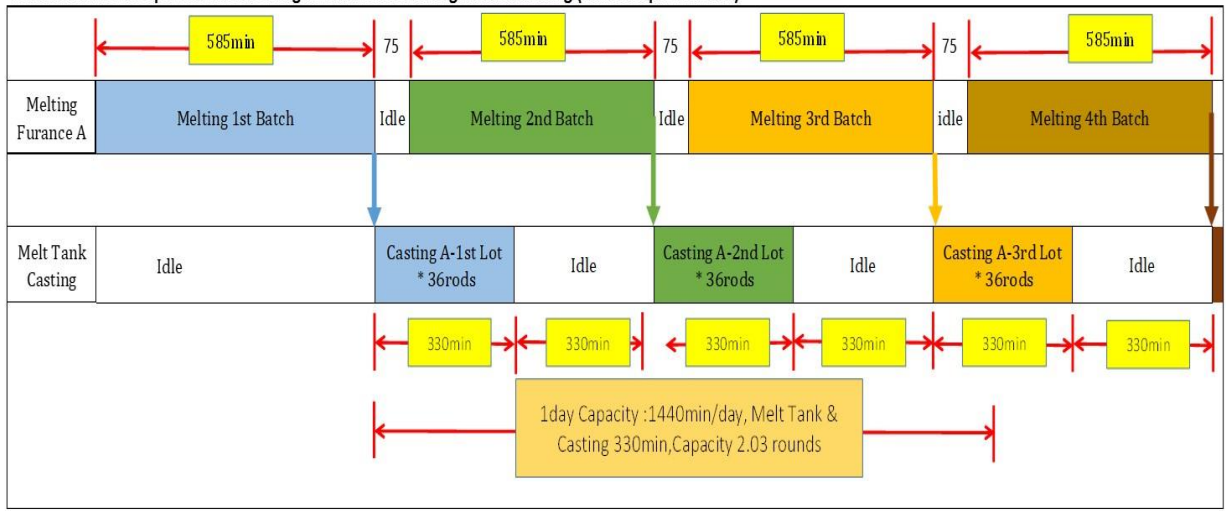

Figure 6: Traditional Time Schedule Melting Furnace and Melting Tank & Casting Chart

In Fig. 7, the Time Schedule to Optimize the Melting Furnace and Melting Tank & Casting Chart is shown the Fig. 7, which must be continuous, with no gaps in between. The materials flow through the layout in a loop shape. The hanger line is required to construct the hanger system and equipment. The system is modernized to set up the control device to move the hanger between workstations and provide the just-in-time information to the manufacturing system. The line balancing for the hanger line can be optimised to increase production efficiency by increasing the throughput time based on increasing the capacity of the bottleneck workstations in the process as the Lean Methodology for Garment Modernization that Prof Dr Ray WM Kong mentioned [2].

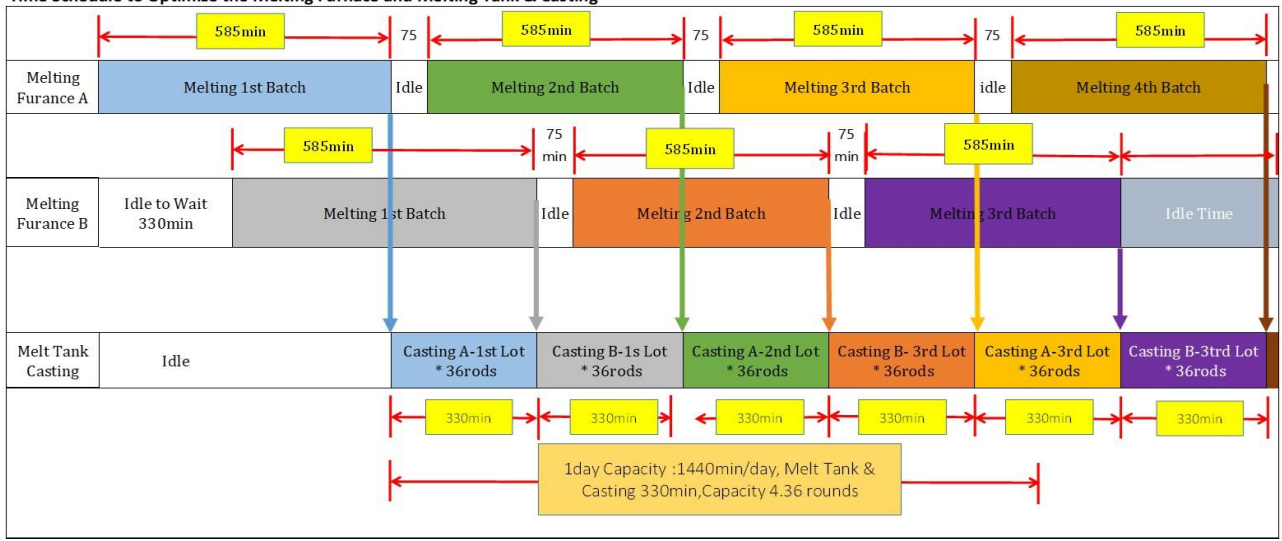

Figure 7: Time Schedule to Optimize the Melting Furnace and Melting Tank & Casting Chart

The plant layout for optimization is required to restructure the melting furnace to the melting tank with the casting facility. The plant layout has added the switch from the melting furnace to the melting tank with the casting facility. Aluminium liquid can be poured from two sets of melting furnaces to one set of melting tanks and a casting facility.

Production output before optimisation: 2 times/day

Production output after optimisation: 4.36 times/day

The increase in output percentage: (4.36 times/day - 2 times/day) / 2 times/day x 100%

The growth output percentage: 118%

There is an optimization in simple line balancing. The labour headcount is increased from 8 workers to 12 workers for 2 shifts a day. The labour cost has increased by 50%. The maximum profit margin of each casting time on 36 rods is over USD28,800 (USD800/rod). The additional profit margin from 2 times to 4.36 times a day is USD67,968. To reduce the additional labour cost for 4 workers, the total amount is USD182/day, so the net profit margin is USD67,786/day.

## VI. LINE BALANCING MODEL ANALYSIS

The main objective of balancing the line is to maximize the marginal profit, which requires reducing the idle cost of melting furnaces. The increased daily output rates of aluminium alloy from the melting tank and the casting facility can provide marginal profit, which also reduces the idle time of melting furnaces. The idle melting furnace is required to provide the energy to keep the temperature of the melting furnace because the drop temperature causes the aluminium liquid to solidify at the bottom of the aluminium residue. The change in temperature causes the melting furnace to crack based on the metallic structure. To keep the temperature from slowing down to drop the temperature, the furnace is not destroyed by thermoshocking, the energy is required to control and keep the expected temperature when the melting furnace is idle.

Referring to the mixed integer linear program (MILP) for the garment line balancing, the optimization of linear programming is to find the optimized solution from Prof Dr Ray WM Kong [3], the number of melting cycles in the melting process can be optimizied to maximise the marginal profit.

Balancing in the stage of a single model finds a locally optimized solution in an iteration. The objective of the stage is to find a solution(s) with a specified number of stations with a minimum cycle time. Solutions are considered locally optimized as the principal objective is to find a solution which will define a smooth production by minimizing the objective function of melting workstation balance from Waldemar Grzechca [8]. The concept of ALBP, where the aim is to optimize the number of workstations with a predefined fixed cycle time, is utilized in the formulation.

Linear programming is adopted to resolve the line balancing with constraints. The formula of linear programming is shown in the following standard format:

$$max_{X_1,X_2,\ldots,X_n} (z) = p_1 x_1 + \ldots + p_n x_n \tag{1}$$

$$\text{subject to} \quad A_{11} x_1 \ldots + A_{1n} x_n \leq b_1,$$

$$\vdots \qquad \qquad \vdots$$

$$A_{m1} x_1 \ldots + A_{mn} x_n \leq b_m,$$

$$x_1, x_2, \ldots, x_n \geq 0$$

By grouping the variables $x_1, x_2, \ldots, x_n$ into a vector $x$ and constructing the following matrix and vectors from the problem data, we can restate the standard form compactly as follows:

$$max_x (z) = p'x$$

$$\text{subject to } Ax \leq b, \ x \geq 0$$

The maximum profit and cycle time is shown in the formula below:

$$Profit_{max} = \max(P_{XY} - \sum_{x=1}^{x=n} IC_x - \sum_{y=1}^{y=n} IC_y) \tag{2}$$

$$Max\ (Profit) = \max(P_{XY} - \sum_{x=1}^{x=n} IC_x - \sum_{y=1}^{y=n} IC_y) \qquad \forall CT_x \in +R, \ \forall CT_y \in +R \tag{3}$$

$$P_{XY} = TC_{XY} * Rr \tag{4}$$

$$TC_{XY} = (CT_x * Eff_x\% + CT_y * Eff_y\%) * Eff_{cast}\% \tag{5}$$

Where $P_{XY}$ is the marginal profit of melting tank and casting in USD,

         $TC_{XY}$ is the output quantity of melting tank and casting in tons,

         $CT_x$ is the cycle time of melting furnace X per cycle in minutes,

         $CT_y$ is the cycle time of melting furnace Y per cycle in minutes,

$Rr$ is the marginal profit rate of melting tank and casting in USD per ton,

$IC_x$ is the idle cost of melting furnace X in minutes,

$IC_y$ is the idle cost of melting furnace Y in minutes,

$Eff_x\%$ is the efficiency of output of melting furnace X after reduction of lose during manufacturing,

$Eff_y\%$ is the efficiency of output of melting furnace Y after reduction of lose during manufacturing,

$Eff_{cast}\%$ is the efficiency of output of melting tank and casting after reduction of lose during manufacturing,

$$\sum_{x=1}^{x=n} IC_x = \sum_{1}^{n}(Cap_x - CT_x\ RT_x)\ CR_x \qquad \forall CT_x \in +R \qquad (6)$$

where $Cap_x$ is the daily available hours of melting furnace X in minutes,

$CT_x$ is the cycle time of melting furnace X per cycle in minutes,

$RT_x$ is the number of time of melting furnace X in minutes per day,

$CR_x$ is the idle cost rate of melting furnace X in USD per minute,

$$\sum_{y=1}^{y=n} IC_y = \sum_{1}^{n}(Cap_y - CT_y\ RT_y)\ CR_y \qquad \forall CT_y \in +R \qquad (7)$$

where $Cap_y$ is the daily available hours of melting furnace Y in minutes,

$CT_y$ is the cycle time of melting furnace Y per cycle in minutes,

$RT_y$ is the number of time of melting furnace Y in minutes per day,

$CR_y$ is the idle cost rate of melting furnace Y in USD per minute,

subject to

$$CT_x\ RT_x \leq\ Cap_x \qquad (8)$$
$$CT_y\ RT_y \leq\ Cap_y \qquad (9)$$

Firstly, the $Profit_{max}$ is required to maximize the daily marginal profit as our goal. The $P_{XY}$ is the profit of aluminium rod per day, which is $P_{XY} = TC_{XY} * Rr$, the $TC_{XY}$ total casting ouput per tank and casting equipment multiple the number of cycle in $Rr$ (number of cycle). The profit margin ($TC_{XY}$) for each casting time means that the casting equipment can produce the 36 rods each cycle, so the $TC_{XY}$ is USD800/rod multiple 36 rods to be USD28,800 (USD800/rod) per cycle of casting.

The $TC_{XY}$ in the formula (5) is the discount the efficiency percentage in the casting facility and equipment, $Eff_{cast}\%$. The melting furnace *x* and the melting furnace *x* and furnace *y* have their efficiency of output of melting furnace X after reduction of loss during the manufacturing process.

The $\sum_{x=1}^{x=n} IC_x$ is the sum idle labor cost of melting furnace *x*. The melting furnace *x* is 1 set in the case study. If the melting furnace *x* is more than 1 set of melting furnaces *x* to n sets, the total idle labor cost is summarized to n sets of melting furnaces *x*. The $\sum_{1}^{n}(Cap_x - CT_x\ RT_x)\ CR_x$ means that the available of daily capacity, $Cap_x$ in the melting furnace *x* minutes the utilized the metling furnace cycle time of melting furnace *x* per cycle in minutes $CT_x$ is multiplied by the number of metling furnace *x*, $RT_x$ and then the result $Cap_x - CT_xRT_x$ is utilizated total idle daily minutes to multiply by the idle cost $CR_x$ for the furnace *x*. The summarized the total idle cost for all furnace *x* in the *n* sets.

The idle cost rate of melting furnace *x*, $CR_x$ is the standard idle labour cost rate in USD per minute when it finds out the idle time in minutes to calculate the idle cost of furnace *x*.

In the constraint, the total daily utilization time is the cycle time of furnace *x* to multiply the number of time of melting furnace *x* in minutes per day, $RT_x$. The daily total utilization time is less than daily capacity of furnace, so it calls the $CT_x \, RT_x \leq Cap_x$.

The $\sum_{y=1}^{y=n} IC_y$ is the sum idle labor cost of melting furnace *y*. The melting furnace *y* is 1 set in the case study. If the melting furnace *y* is more than 1 set of melting furnaces *y* to n sets, the total idle labor cost is summarized to n sets of melting furnaces *y*. The $\sum_1^n (Cap_y - CT_y \, RT_y) \, CR_y$ means that the available of daily capacity, $Cap_y$ in the melting furnace *y* minutes the utilized the melting furnace cycle time of melting furnace *y* per cycle in minutes $CT_y$ is multiplied by the number of metling furnace *y*, $RT_y$ and then the result $Cap_y - CT_y RT_y$ is utilizated total idle daily minutes to multiply by the idle cost $CR_y$ for the furnace *y*. The summarized the total idle cost for all furnace *y* in the *n* sets.

The idle cost rate of melting furnace *y*, $CR_y$ is the standard idle labour cost rate in USD per minute when it finds out the idle time in minutes to calculate the idle cost of furnace *y*.

In the constraint, the total daily utilization time is the cycle time of furnace *y* multiplied by the number of times of melting furnace *y* in minutes per day, $RT_y$. The daily total utilization time is less than daily capacity of furnace, so it calls the $CT_y \, RT_y \leq Cap_y$.

Hence, the above formula from (1) to (9) has shown the detail of calculation for the linear progamming in the line balancing.

In the case study, the increase in daily marginal profit is attributed to the enhanced output from melting and casting processes, as demonstrated by the following calculation:

The increase in daily marginal profit is calculated as:

> Increment in daily marginal profit = (New optimized daily cycles - Original daily cycles) *Output cost per cycle per day - Additional operators' cost

Substituting the values:

> = (4.36 cycles – 2 cycles) * USD28,800/cycle/day – (USD182/day)
>
> = USD67,968/day – USD182/day
>
> = USD67,786/day

The growth in marginal daily profit is calculated as:

> = (USD125,386 – USD57,600) / USD57,600
>
> = 117.6%

In this optimization scenario, simple line balancing was achieved by increasing the labour headcount from 8 to 12 workers across two shifts per day, resulting in a 50% increase in labour costs. The maximum profit margin per casting cycle of 36 rods exceeds USD28,800 (USD800 per rod). The additional profit margin from increasing production from 2 to 4.36 cycles per day amounts to USD67,968. After accounting for the additional labour cost of USD182 per day for the 4 extra workers, the net profit margin is USD67,786 per day, reflecting a 117.6% growth in marginal daily profit.

## VII. CONCLUSION

Linear programming, particularly mixed integer linear programming (MILP), offers several benefits for line balancing in aluminium extrusion processes. Here are the key benefits and conclusions of applying linear programming to this research and development of the aluminium extrusion industry:

1. Optimization of Resources:

- Linear programming helps in optimizing the allocation of resources, such as melting furnaces and casting facilities, to maximize output and minimize idle time. This ensures that resources are used efficiently, reducing waste and operational costs.

2. Maximization of Marginal Profit:

- By optimizing the number of melting cycles and balancing the line, linear programming can help increase daily output rates and greater utilization, thereby calculating the maximizing marginal profit. This is achieved by reducing idle labour costs and ensuring continuous operation of melting furnaces $x$ and $y$.

3. Reduction of Idle Time:

- The approach minimizes idle time for melting furnaces (furnaces $x$ and $y$ in n furnace), which is crucial for maintaining the necessary temperature to prevent aluminium solidification and potential damage to the furnace as a crash of furnace and melting tanks. This leads to more consistent production and less downtime.

4. Improved Production Flow:

- Linear programming aids in achieving a smooth production flow by determining the optimal number of cycle times of the melting process. This reduces bottlenecks and ensures a balanced workload across the production line.

5. Enhanced Decision-Making:

- The use of linear programming provides a data-driven approach to decision-making, allowing aluminium extrusion managers to evaluate different scenarios and choose the most efficient production strategy.

6. Flexibility and Scalability:

- The methodology can be adapted to various production models, whether single-model, mixed-model, or multi-model lines, providing flexibility to accommodate changes in production demands or product types.

The application of linear programming, particularly MILP, in aluminium extrusion line balancing, offers significant advantages in optimizing production processes. By focusing on maximizing marginal profit and minimizing idle costs, linear programming ensures efficient resource utilization and enhances overall productivity. The approach provides a structured framework for balancing production lines, reducing idle time, and maintaining the operational integrity of melting furnaces. Ultimately, this leads to improved profitability, reduced operational risks, and a more agile production system capable of adapting to market demands and technological advancements. In conclusion, effective line balancing in garment assembly operations is essential for resolving the issue of excessive work-in-process (WIP) inventory and improving the production output and efficiency, which often arises from an unbalanced production line. By systematically analyzing and optimizing the distribution of tasks among operators, organizations and factory planners can achieve a more synchronized workflow that minimizes bottlenecks and reduces idle time. This not only leads to a smoother production process but also significantly decreases the accumulation of WIP, thereby lowering storage costs and enhancing overall operational efficiency. Ultimately, implementing line-balancing techniques fosters a leaner manufacturing environment, improves responsiveness to market demands, and contributes to higher levels of productivity and profitability in the garment industry.

The research article on aluminium extrusion can draw insights from cross-industry studies, highlighting the critical role of automation and intelligent manufacturing in boosting productivity and output. The development of innovative technologies, such as the "Design and Experimental Study of Vacuum Suction Grabbing Technology to Grasp Fabric Piece" by Prof. Dr. Ray WM Kong et al. [9] and in the K. M. Batoo (Ed.), Science and Technology: Developments and Applications, Innovative Vacuum Suction-grabbing Technology for Garment Automation from Prof. Dr. Ray WM Kong et al. [10]. Similarly, the "Design of a New Pulling Gear for the Automated Pant Bottom Hem Sewing Machine" by Prof. Dr. Ray WM Kong et al. [11] demonstrates how automation can enhance production rates in hem sewing machines. Looking ahead, automation presents a significant opportunity to advance the aluminium extrusion industry, paralleling its transformative impact on garment and electronics manufacturing sectors.